# Taylor Series as Wide-sense Biorthogonal Wavelet Decomposition

H.M. de Oliveira, R.D. Lins

*Abstract*—**Pointwise-supported generalised wavelets are introduced, based on Dirac, doublet and further derivatives of $\delta(t)$. A generalised biorthogonal analysis leads to standard Taylor series and new Dual-Taylor series that may be interpreted as Laurent Schwartz distributions. A Parseval-like identity is also derived for Taylor series, showing that Taylor series support an energy theorem. New representations for signals called '*derivagrams*' are introduced, which are similar to 'spectrograms'. This approach corroborates the impact of wavelets in modern signal analysis.**

*Index Terms*—**Pointwise-supported wavelets, wide-sense biorthogonal analysis, Parseval theorems, derivagrams.**

## I. INTRODUCTION

Continuous and discrete wavelet transforms have emerged as a widely used tool in signal analysis, which has been proved to be valuable for scientists and engineers [1-2]. Today, they are applied extensively in applications in an amazing number of areas [3-5], including seismic geology, image processing (e.g., video data compression, 'denoising'), target recognition and radar, spectrometry, metallurgy, turbulence, computer graphics, transient analysis in power lines, computer and human vision, optics and electromagnetism, telecommunications, DNA sequence analysis, fractals, acoustic signal characterisation, quantum physics, biomedical signal analysis (mammography, electrocardiogram, etc.), earthquake forecast, statistics, solution of partial and ordinary differential equations, amongst many others. Orthogonality has long been assumed as a key property in virtually all standard approaches when analysing or synthesising signals [1,3]. A higher-level signal processing technique involves the concept of biorthogonality in which two (cross-orthogonal) sets are used – one in the analysis and the other in the synthesis of signals [2]. In the late 1990's, biorthogonal Wavelets brought a major breakthrough into image compression, thanks to their natural feature of concentrating energy in a few transform coefficients. It has been adopted in international standards such as JPEG 2000 [6,7] and the FBI standard for fingerprint storage [8,9]. A recent comparative analysis of the file formats and algorithms for still image compression shows the unbeatable superiority of JPEG 2000 [23] for image storage and network transmission. The aim of this paper is to investigate new wide-sense biorthogonal decompositions.

Federal University of Pernambuco, Signal Processing Group, C.P. 7.800, 50.711-970, Recife – PE. This work was partially supported by the Brazilian National Council for Scientific and Technological Development (CNPq) under research grants #.306180 (HMdO) and #.523974/96-5 (RDL), e-mail: {hmo,rdl}@ufpe.br

The paper is organised as follows. This Section briefly overviews wavelets and some relevant aspects of the Laurent Schwartz theory of distributions [10-12], as well as their relationship with wavelets. Section II addresses standard Taylor series viewed as Tukey-tapers [13]. The new wide-sense wavelets derived from Dirac derivatives are introduced in Section III, showing that Taylor series can be interpreted as some sort of biorthogonal analysis based on impulsive wavelets. A Parseval-like identity for infinite Taylor series is then derived and a few examples are presented. Section IV introduces the Taylor-energy density and a new signal analysis tool called '*derivagram*'. Final remarks are presented in the last Section, which emphasizes the power and applicability of wavelet theory in modern signal analysis.

Here, the symbol ":=" denotes "equal by definition". $\mathbb{N}, \mathbb{R}$ and $\mathbb{C}$ are the set of Natural, Real and Complex numbers, respectively. Wavelets are denoted by $\psi(t)$ and scaling functions by $\phi(t)$, with corresponding spectra $\Phi(w)$ and $\Psi(w)$, respectively.

**Definition 1** ($C^\infty$ space). The space of infinitely differentiable complex signals $f: \mathbb{R} \to \mathbb{C}$ is denoted by $C^\infty(\mathbb{R})$, i.e. the vector space such that
$$\forall f \in C^\infty \Rightarrow (\forall n \in \mathbb{N}) f^{(n)} \in C^\infty,$$
where $f^{(n)}$ denotes the $n^{th}$ derivative of $f$. ∎

**Definition 2** (Support of a signal). The support of a given signal $f: \mathbb{R} \to \mathbb{C}$, supp $f$, is the closed union of the set
$t \in \mathbb{R} \mid f(t) \neq 0$. Rigorously, supp $f := \overline{\{t \in \mathbb{R} \mid f(t) \neq 0\}}$. ∎

The set of interest for distributions is the set D.

**Definition 3** (D space). Let D be the vector space of signals $\varphi: \mathbb{R} \to \mathbb{C}$ with infinitely derivatives of bounded support. ∎

A distribution is a continuous linear functional over D, defined by an assignment rule
$$T : D \to C$$
$$\varphi \mapsto T(\varphi).$$

Notation: Let D' be the space of all distributions.

*Dirac Distribution* ($\delta$). $\forall \varphi \in D$,
$$<\varphi, \delta> := \int_{-\infty}^{+\infty} \varphi(t) \delta(t) dt := \varphi(0). \blacksquare$$

A Dirac distribution evaluated on a point $t_0 \in \mathbb{R}$ is defined by
$$\forall \varphi \in D, \quad <\varphi, \delta_{(t_0)}> := \int_{-\infty}^{+\infty} \varphi(t) \delta_{(t_0)}(t) dt := \varphi(t_0).$$

Usually $T(\varphi) = <\varphi, T>$ is defined via an inner product, thus being a linear operator.

*Derivative of Distributions.* $<\varphi, T'> := -<\varphi', T>$.

Higher-order derivatives (order *n*). By recurrence, the $k^{th}$ derivative ($k>0$) of $T \in D'$, denoted by $T^{(k)}$, is given by:

$$\forall \varphi \in D \quad <\varphi, T^{(k)}> = (-1)^k <\varphi^{(k)}, T>. \blacksquare$$

This definition supply the following corollary:
**Corollary 1**. Every distribution is infinitely derivable.

**Definition 4** (The support of a distribution). The support of a distribution $T$, denoted supp $T$, is the smallest closed set in which $T$ is not vanishing, that is, $\exists \varphi \in D \mid <\varphi, T> \neq 0$. $\blacksquare$

*Compacted supported distributions*. Let $T$ be a distribution of bounded support, thus compact. One can prove [10] that $T$ can be extended to the space $D:=C^\infty(\mathbb{R})$. The set of compactly supported distributions is denoted by $D'$. $\blacksquare$

Continuous wavelet transforms [1,4] can also be interpreted as distributions [10,11]. If $\psi(t)$ is a mother wavelet, then:

$$\forall \varphi \in D, \quad <\varphi, CWT> = \int_{-\infty}^{+\infty} \varphi(t) \psi(t) dt.$$

The homothetic and translation relationships of distributions [10] correspond to usual scaling and translation in wavelets.

$$i) <\varphi, \tau_b \psi> = \int_{-\infty}^{+\infty} \varphi(t) \psi(t-b) dt.$$

$$ii) \; a \neq 0 \quad <\varphi, \psi_a> = \int_{-\infty}^{+\infty} \varphi(t) \frac{1}{|a|} \psi(\frac{t}{a}) dt.$$

Therefore, in the general case, a wavelet transform corresponds to a distribution:

$$<\varphi, \tau_b \psi_a> = \int_{-\infty}^{+\infty} \varphi(t) \frac{1}{|a|} \psi(\frac{t-b}{a}) dt.$$

## II. Taylor Data Window

The standard Taylor series [14], truncated with $N+1$ terms, can be interpreted as a linear functional over the function. If $f$ possesses $N$ derivatives at the origin, then let us define the Taylor kernel as,

$$K_N(t;t') := \sum_{n=0}^{N} \frac{(-t)^n}{n!} \delta^{(n)}(t').$$

Clearly, $\int_{-\infty}^{+\infty} K_N(t;t') dt' = 1$. The truncated Taylor series $f_N(t)$ can be obtained by the observation of $f(t)$ thought the Taylor *taper* [13]:

$$f_N(t) := \int_{-\infty}^{+\infty} f(t') K_N(t;t') dt'.$$

The use of tapers is currently being investigated to devise new multiplex systems (private communication). It is easy to show that

$$f_N(t) = \sum_{n=0}^{N} \frac{f^{(n)}(0)}{n!} t^n$$

such that the limit,

$$\lim_{N \to \infty} f_N(t) = f(t)$$

provided that the Taylor series converges. The set of values of the independent variable *t*, for which the series converge, constitutes what is called the region of convergence of the series.

## III. Generalised Impulsive Wavelets

**Definition 5** (Generalised wavelets of Dirac). We define the "impulsive wavelets" at a scale $a \in \mathbb{N}\text{-}\{0\}$ and a translation $b \in \mathbb{R}$ as,

$$\psi_{a,b}(t) := (-1)^a \delta^{(a)}(t-b), \quad -\infty < t < +\infty \quad \blacksquare$$

The scale parameter *a* plays here a role somewhat different of the one it stands for in conventional wavelet theory.
These "wavelets" have compact support and their narrow support is given by Supp $\psi_{a,b}(t) = \{b\}$. Thus, they are pointwise-supported generalised wavelets. Let $\psi_{a,b}(t) \leftrightarrow \Psi_{a,b}(w)$ be a pair formed by a wavelet and its Fourier transform. Since their generalised spectra are $\Psi_{a,b}(w) = (-jw)^a e^{-jwb}$, these are neither signals of finite energy nor hold the admissibility condition [1,2]. Nevertheless, one can remark that

$$\int_{-\infty}^{+\infty} \psi_{a,b}(t) dt = 0, \; \forall \; a \in \mathbb{N}\text{-}\{0\}, b \in \mathbb{R}.$$

These generalized-wavelets are associated with a Dirac scale function $\phi(t) = \delta(t)$. Clearly, $\int_{-\infty}^{+\infty} \phi(t) dt = 1$, as expected. An "extensive" multiresolution analysis [2,15] can as a result be implemented using the following analysing sets:

$$\phi = \{\delta(t)\}; \text{ generalised scale function,}$$

$$\psi = \{-\delta'(t), \delta''(t), -\delta'''(t), \ldots\}; \text{ generalised wavelets at different scales.}$$

If *f(t)* has derivatives of all orders throughout a neighbourhood of a point $t_0=b$, then for any continuous signal $f(t) \in C^\infty$,
the wavelet coefficients are given by

$$\forall \; a \in \mathbb{N}\text{-}\{0\}, b \in \mathbb{R}$$

$$c_{a,b} = WT_{a,b}[f(t)] := <f(t), (-1)^a \delta^{(a)}(t-b)>$$

$$= \int_{-\infty}^{+\infty} f(t)(-1)^a \delta^{(a)}(t-b) dt$$

$$= f^{(a)}(b),$$

that is, they correspond to the $a^{th}$ derivative of the analysing signal *f(t)* at the point *t=b*.

The reconstruction guidelines can be derived through the following biorthogonal property of the set $\{\psi\}$:

Let the reconstruction functions be defined by:

$$\widetilde{\phi} = \{\widetilde{\phi}(t) := 1\}; \; \widetilde{\psi} = \left\{\widetilde{\psi}_{a,b}(t) := \frac{(t-b)^a}{\Gamma(a+1)}\right\}.$$

Since *a* is constrained to be an integer, the set of interest is

$$\left\{1, \frac{t-b}{1!}, \frac{(t-b)^2}{2!}, \frac{(t-b)^3}{3!}, \ldots\right\}.$$

These reconstruction functions are not wavelets, but hold the biorthogonality property.

**Proposition 1** (biorthogonality). $\forall n,m \in \mathbb{N}\text{-}\{0\}, t_0 \in \mathbb{R}$

$$\int_{-\infty}^{+\infty} \psi_{n,t_0}(t) \cdot \widetilde{\psi}_{m,t_0}(t) dt = \delta_{n,m},$$

where $\delta_{n,m}$ is the Kronecker symbol.

Proof:

$$<\psi_{n,t_0},\widetilde{\psi}_{n,t_0}> = \int_{-\infty}^{+\infty} \frac{(t-t_0)^m}{m!} \cdot (-1)^n \delta^{(n)}(t-t_0)dt$$

$$= \frac{1}{m!} \cdot \frac{d^n}{dt^n}\left((t-t_0)^m\right)\bigg|_{t=t_0} = \delta_{n,m}. \blacksquare$$

The generalised reconstruction formula is therefore

$$f(t) \cong c_{0,t_0}\widetilde{\phi}(t-t_0) + \sum_{n=1}^{+\infty} c_{n,t_0}\widetilde{\psi}_{n,t_0}(t),$$

which yields,

$$f(t) = f(t_0) + \sum_{n=1}^{+\infty} \frac{f^{(n)}(t_0)}{n!} \cdot (t-t_0)^n,$$

a standard Taylor series at a particular point $t_0$ that corresponds exactly to the $b$ parameter of the analysis. The equality holds within the interval of convergence of the Taylor series [14]. Besides consolidating wavelet representation as a most powerful tool, this approach suggests the development of a novel Taylor-like series, based on further biorthogonal analysis [15].

We may also set

$$\widetilde{c}_{a,b} := <f(t),\widetilde{\psi}_{a,b}(t)> = \int_{-\infty}^{+\infty} f(t) \cdot \frac{(t-b)^a}{a!} dt.$$

Therefore $f$ can be approximated by

$$f(t) \cong c_{0,t_0}\phi(t-t_0) + \sum_{n=1}^{+\infty} c_{n,t_0} \cdot \widetilde{\psi}_{n,t_0}(t).$$

For the sake of simplicity, we denote

$$f^{(\widetilde{n})}(t_0) := \int_{-\infty}^{+\infty} f(t) \cdot (t-t_0)^n dt.$$

When $f(t)$ is a probability distribution, this parameter corresponds to the $n^{\text{th}}$ moment of the density $f(t)$ [14].

The dual of the Taylor series is the inhomogeneous representation [2]

$$f(t) \stackrel{!}{=} f^{(\widetilde{0})}(t_0)\delta(t-t_0) + \sum_{n=1}^{+\infty} \frac{f^{(\widetilde{n})}(t_0)}{n!} \cdot (-1)^n \delta^{(n)}(t-t_0).$$

The series above is quite unconventional and the equality "$\stackrel{!}{=}$" must be suitably interpreted. This is an extensional concept as identity of distributions. Under integration, the term at the right hand side provides the same result as the integration of $f(t)$ itself. By no means it implies that the two members of the equation are identical (pointwise convergence).

A new Parsevel-Taylor identity can be finally established (proof omitted here):

**Proposition 2** (Parsevel-Taylor). If $f \in L^2$ is a real signal and $f \in C^\infty$, then

$$\int_{-\infty}^{+\infty} f^2(t)dt = \sum_{n=0}^{+\infty} c_{n,b} \cdot \widetilde{c}_{n,b} = \sum_{n=0}^{+\infty} \frac{f^{(n)}(b) \cdot f^{(\widetilde{n})}(b)}{n!}. \blacksquare$$

$E := \int_{-\infty}^{+\infty} f^2(t)dt$ is the signal energy.

**Corollary 2**. If the region of convergence of the abovementioned series is RC $\subset$ (-∞,+∞), then

$$\int_{RC} f^2(t)dt = \sum_{n=0}^{+\infty} c_{n,b} \cdot \widetilde{c}_{n,b} = \sum_{n=0}^{+\infty} \frac{f^{(n)}(b) \cdot f^{(\widetilde{n})}(b)}{n!},$$

where simply the moments are now

$$f^{(\widetilde{n})}(t_0) := \int_{RC} f(t) \cdot (t-t_0)^n dt. \blacksquare$$

Energy theorems are among the chief results of series decomposition and finite or infinite transforms [2,4-5]. Standard Taylor series has not got, to the best of our knowledge, a similar relationship.

**EXAMPLE 1.** Let $f(t) = e^{-t^2} \in C^\infty(\mathbb{R})$ be the Gaussian pulse. Indeed, $f \in L^2(\mathbb{R})$, as

$$\int_{-\infty}^{+\infty} f^2(t)dt = \sqrt{\pi}.$$

Evaluating the Taylor series at the origin, one obtains

$$f(t) = 1 - \frac{t^2}{2} + \frac{t^4}{2^2 \cdot 2!} - \frac{t^6}{2^3 \cdot 3!} + \ldots + R_n.$$

The derivatives of $f(t)$ can easily be derived from the series above, $\forall k \in \mathbb{N}$,

$$\begin{cases} f^{(2k+1)}(0) = 0 \\ f^{(2k)}(0) = \frac{(-1)^k (2k)!}{2^k \cdot k!} \end{cases}$$

Let $\Gamma(.)$ be the Factorial (gamma) function or Euler's integral of the second kind [16],

$$\Gamma(x) := \int_0^{+\infty} e^{-t} t^{x-1} dt, \; t>0.$$

The moments of the Gaussian pulse are given by [14,16], $\forall k \in \mathbb{N}$,

$$\begin{cases} f^{(\widetilde{2k+1})}(0) = 0 \\ f^{(\widetilde{2k})}(0) = \int_{-\infty}^{+\infty} t^{2k} e^{-t^2/2} dt = 2^{k+\frac{1}{2}} \Gamma\left(k + \frac{1}{2}\right) \end{cases}$$

The dual Taylor series for the Gaussian pulse is therefore

$$e^{-t^2/2} \stackrel{!}{=} \sqrt{2\pi}\delta(t) + \sum_{k \text{ even}} 2^{k+\frac{1}{2}} \Gamma\left(k + \frac{1}{2}\right) \delta^{(k)}(t).$$

The Parseval-Taylor identity can be applied to derive the following non-trivial numerical series:

$$\sum_{k \text{ even}}^{+\infty} (-1)^k \frac{\Gamma\left(k + \frac{1}{2}\right)}{\Gamma(k+1)} = \sqrt{\frac{\pi}{2}}.$$

The convergence of the series is a bit slow as shown in Fig.1. Let

$$\text{converg1}(N) := \sum_{k=0}^{N} (-1)^k \cdot \frac{\Gamma(k + 0.5)}{\Gamma(k + 1)}$$

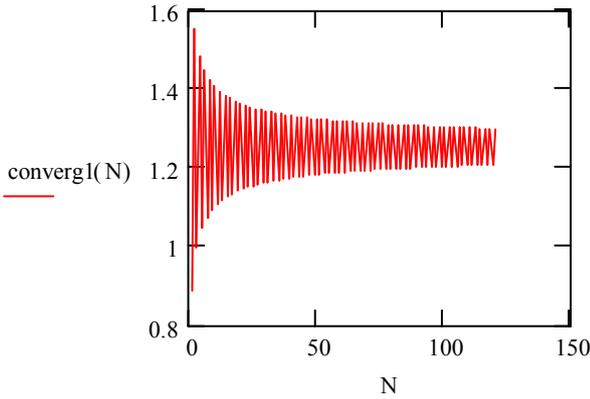

**Figure 1**. The convergence of the numerical series corresponding to the Parseval-Taylor identity applied to the Gaussian pulse.
$\Gamma(.)$ is the generalized factorial function. The series converges (oscillatory behaviour) to $\sqrt{\pi/2} = 1.25331414...$

If $f(t) \leftrightarrow F(w)$ denotes the Fourier pair, then
$$f^{(\tilde{n})}(0) = \int_{-\infty}^{+\infty} t^m f(t)dt = Fourier\left(t^m f(t)\right)\Big|_{w=0} = j^n F^{(n)}(0).$$
As a result, both terms $f^{(n)}(0)$ and $f^{(\tilde{n})}(0)$ of the analysis are related to derivatives of the signal, in the time and in the frequency domain, respectively.

**Proposition 3** (moment symmetry). Even (respectively odd) signals have only non-zero even (respectively odd) moments, i.e. If $f \in C^\infty(-\infty,\infty)$ is even, then
$$\forall k, \; f^{(\widetilde{2k+1})}(0) = F^{(2k+1)}(0) = 0.$$
Conversely, If $f \in C^\infty(-\infty,\infty)$ is odd, then
$$\forall k, \; f^{(\widetilde{2k})}(0) = F^{(2k)}(0) = 0.$$
Proof. $f^{(\tilde{n})}(0) = \int_{-\infty}^{+\infty} t^m f(t)dt = (-1)^m \int_{-\infty}^{+\infty} t^m f(-t)dt$ and the proof follows. ∎

**Corollary 3**. The above proposition also holds if
$$\exists L < +\infty \; | \; f \in C^\infty(-L,L). \; \blacksquare$$

**Corollary 4.** The Taylor-Parseval identity (proposition 2) can be rewritten as
$$\sum_{n=0}^{+\infty} j^n \cdot \frac{f^{(n)}(0)F^{(n)}(0)}{n!} = \int_{-\infty}^{+\infty} f^2(t)dt. \; \blacksquare$$
Symmetric energy propositions.
$$\sum_{n=0}^{+\infty} (-1)^n \cdot \frac{f^{(2n)}(0)F^{(2n)}(0)}{2n!} = \int_{-\infty}^{+\infty} (E(f(t)))^2 dt$$
and
$$\sum_{n=0}^{+\infty} (-1)^n \cdot j \cdot \frac{f^{(2n+1)}(0)F^{(2n+1)}(0)}{(2n+1)!} = \int_{-\infty}^{+\infty} (O(f(t)))^2 dt.$$

**EXAMPLE 2.** Let now $f(t) = \ln(1-t) \in C^\infty(-1,1)$. The well-known Taylor series is [14,16]
$$\ln(1-t) = -\left(t + \frac{1}{2}t^2 + \frac{1}{3}t^3 + \frac{1}{4}t^4 + ...\right), \quad -1 \leq t < 1.$$
Now, $f(t) \in L^2(-1,1)$ since
$$\int_{-1}^{1} [\ln(1-t)]^2 dt = 2\ln 2\{\ln 2 - 2\} + 4.$$
Indeed, $f(t) \in L^2(0,1)$ and $\int_0^1 [\ln(1-t)]^2 dt = 2$.

The derivatives of $f(t)$ at the origin are given by $\forall n \in \mathbb{N}-\{0\}$
$$\begin{cases} f^{(n)}(0) = -\dfrac{\Gamma(n+1)}{n} \\ f^{(0)}(0) = 0. \end{cases}$$

Let us examine a sub-region of convergence $(0,1) \subset (-1,1)$ as a particular case. The moments of $f(t)$ can be computed in terms of Euler's psi function (logarithmic factorial or digamma function [16])
$$\psi(x) := \frac{d}{dx} \ln \Gamma(x).$$
The value $C := -\psi(1) = 0.57721566490...$ is the Euler's constant. Thus,
$$f^{(\tilde{n})}(0) = \int_0^1 t^n \ln(1-t) dt = -\frac{1}{n+1}[\psi(n+2) - \psi(1)].$$
A rational function related to the harmonic function is $F(x) := \psi(x+1) + C$.
If $x$ is an integer, it follows then that
$$F(n) = \psi(n+1) + C = \sum_{k=1}^{n} \frac{1}{k},$$
so $F$ is closed-linked to the harmonic series. The moments of this logarithmic pulse can be written as
$$f^{(\tilde{n})}(0) = -\frac{F(n+1)}{n+1},$$
and the Parseval-Taylor series yields the following non-trivial identity:
$$\sum_{n=1}^{+\infty} \frac{F(n+1)}{n(n+1)} = 2 \; \text{or} \; \sum_{n=1}^{+\infty} \frac{\psi(n+2)}{n(n+2)} = 2 - C,$$
$$\sum_{n=1}^{+\infty} \frac{(-1)^n F(n+1)}{n(n+1)} = 2(\ln^2 2 - 2\ln 2 + 1).$$

The convergence of one amongst the series above is illustrated in Fig. 2.

Let
$$\text{converg2}(N) := \sum_{k=1}^{N} \left[\frac{\text{Psi}(k+2)}{k \cdot (k+1)}\right].$$

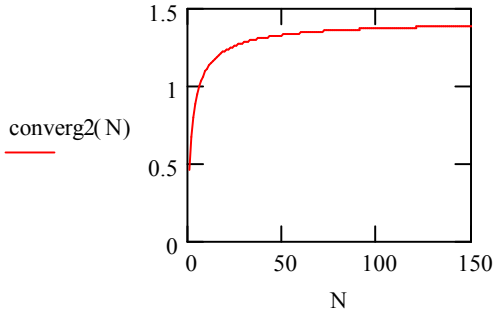

**Figure 2**. The convergence of the numerical series corresponding to the Parseval-Taylor identity applied to the logarithmic pulse. Psi(.) is the Euler's digamma function. The series monotonically converges to 2-*C*=1.42278434…, where *C* is the Euler's constant.

$$\int_{-1}^{1}\left(\frac{1}{2}\ln(1-t^2)\right)^2 dt + \int_{-1}^{1}\left(\frac{1}{2}\ln\left(\frac{1-t}{1+t}\right)\right)^2 dt =$$
$$2\ln^2(2) - 4\ln(2) + 4.$$

## IV. 'DERIVAGRAMS' AND TAYLOR-ENERGY DENSITY

New tools for analysing signals can be derived under the present framework, namely the Taylor energy density of a signal and '*derivagrams*'. From Proposition 2, signals of finite energy hold

$$E = \sum_{n=0}^{+\infty} \frac{f^{(n)}(b).f^{(\tilde{n})}(b)}{n!} < +\infty.$$

Consequently, the contribution at the $n^{th}$-derivative level to the full energy of the signal can be defined according to

$$DE_n := \frac{f^{(n)}(b).f^{(\tilde{n})}(b)}{n!}, \quad n \in \mathbb{N}.$$

This quantity will be referred to as the (pseudo) Taylor energy density. Finally, we have the following energy decomposition:

$$E = \sum_{n=0}^{+\infty} DE_n.$$

This is not precisely a density of energy because $DE_n$ may be negative. Since that $f \in L^2(\mathbb{R}C)$, the series above is convergent so that $(\forall n), |DE_n| < +\infty$.

The cumulative energy until the $N^{th}$-derivative level of the signal is defined by

$$E_N := \sum_{n=0}^{N} DE_n.$$

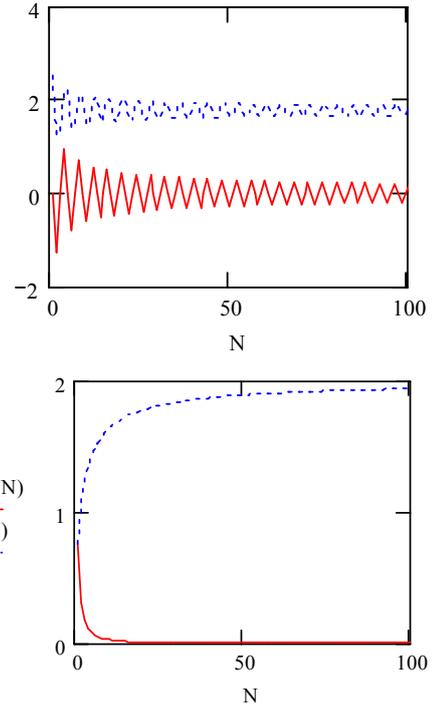

**Figure 3.** Taylor energy density of: (a) gaussian pulse RC $=(\infty,+\infty)$ and (b) logarithmic pulse RC(-1,1). *DE*(*N*) denotes the Taylor energy density at the level *N*, and *E*(*N*) is the cumulated energy due to signal components until the $N^{th}$ derivative.

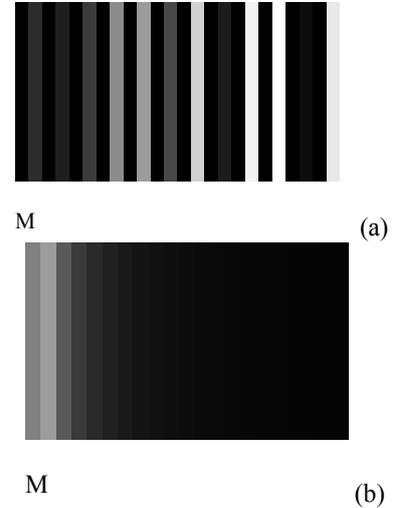

**Figure 4.** 'Derivagrams' with 10 levels (*n*=0,1,..,9) at $t_0$=0 for: (a) the gaussian pulse RC=$(-\infty,+\infty)$ and (b) the logarithmic pulse RC(-1,1). The colour map adopted has 1024-Gray-level. *DE*(*N*) denotes the Taylor energy density at the level *N*, and *E*(*N*) is the accumulated energy due to signal components until the $N^{th}$ derivative. The total energies are $\sqrt{\pi} = 1.77245385...$ and 2, respectively.

**Definition 6** (derivagrams). We can create new pictorial energy representations (in the same way as spectrograms, scalograms) referred to as 'derivagrams', according to the numerical values of

$$derivagram(n, t_0) := \frac{f^{(n)}(t_0).f^{(\tilde{n})}(t_0)}{n!}, n \in \mathbb{N}, t_0 \in \mathbb{R},$$

that can be depicted as a bargraph or a piecewise curve. ∎

Figures 3 and 4 show the 'derivagrams' for the pulses discussed in examples 1 and 2.

## V. Closing Remarks

This paper introduces an original reading for the classical Taylor series, which is based on the wavelet theory. This approach corroborates the power of wavelet analysis, by interpreting Taylor series in the wavelet framework. The wavelet analysis already encompassed tools as Fourier, Gabor, 'spectrograms', 'scalograms' [17], pyramidal algorithm [18], Heisenberg inequalities for Wavelet [19-20], and even Gibbs phenomenon have been addressed [21-22].

A dual version of the Taylor series based on distributions is also presented. A Parseval-Taylor energy theorem is also proposed for Taylor series. These new representations may be valuable for describing features embedded in signals. New signal processing tools were introduced, which support energy theorems based on Taylor series. The present analysis can derive the value of many infinite series in the same way as Fourier series does.

The results presented go far beyond simple mathematical curiosity. They cast some light and generality on the mechanism of signal analysis. Similarly to the widespread idea of vanishing moments often adopted in the wavelet design [1-5], particular series for which $f^{(n)}(0) = 0$ $n=1,2...N$ (Butterworth) should carefully be investigated.

## VI. Appendix

Surprising, the celebrated Shannon/Nyquist/Kotel'nikov sampling theorem [24-25] can also be viewed as a generalised biorthogonal reconstruction.

**Proposition 4.** The sets $\{\psi_n(t)\} := \{Sa(2\pi Bt - n\pi)\}_{-\infty}^{+\infty}$ e $\{\tilde{\psi}_n(t)\} := \left\{\delta\left(t - \frac{n}{2B}\right)\right\}_{-\infty}^{+\infty}$ are biorthogonals. ∎

Proof. One has

$$\int_{-\infty}^{+\infty} \psi_n(t)\tilde{\psi}_m(t)dt = \int_{-\infty}^{+\infty} Sa(2\pi Bt - n\pi)\delta(t - \frac{m}{2B})dt = Sa(\pi(n-m))$$

Therefore, there are the following reconstruction formulas:

i) $f(t) = \sum_{m=-\infty}^{+\infty} \tilde{c}_m \psi_m(t) = \sum_{m=-\infty}^{+\infty} f\left(\frac{m}{2B}\right).Sa(2\pi Bt - m\pi)$,

where $\tilde{c}_m = \int_{-\infty}^{+\infty} f(t)\delta(t - \frac{m}{2B})dt = f\left(\frac{m}{2B}\right)$,

ii) $f(t) = \sum_{n=-\infty}^{+\infty} c_n \tilde{\psi}_n(t) = \sum_{n=-\infty}^{+\infty} \tilde{f}\left(\frac{n}{2B}\right)\delta\left(t - \frac{n}{2B}\right)$,

where $c_n = \int_{-\infty}^{+\infty} f(t).Sa(2\pi Bt - n\pi)dt := \tilde{f}\left(\frac{n}{2B}\right)$.